\documentclass[11pt,a4paper]{article}
\usepackage[T1]{fontenc}
\usepackage[utf8]{inputenc}
\usepackage{amsmath,amssymb}
\usepackage[margin=2.5cm]{geometry}
\usepackage{hyperref}
\hypersetup{
    colorlinks=true,
    linkcolor=blue,
    urlcolor=blue,
    citecolor=blue
}
\usepackage{enumitem}

\title{\textbf{Structure-Driven Methodology: An Emerging Cross-Domain Paradigm}}
\author{CHEN Shengchang\\
 School of Earth Sciences, Zhejiang University, Hangzhou, China\\
\texttt{chenshengc@zju.edu.cn}}
\date{}

\begin{document}

\maketitle

\begin{abstract}
This article proposes and articulates structure-driven methodology as an emerging cross-domain paradigm. In contrast to the objective-driven paradigm that has dominated inverse-problem solving for over two centuries, structure-driven methodology advocates identifying and exploiting the intrinsic structure of a problem to construct the solution method directly, rather than approaching the solution through iterative search. Drawing on examples from the Liouville-variant Goldbach conjecture, the three-dimensional Kakeya conjecture, Noether's theorem, and the author's nearly four decades of geophysical inversion practice, the article illustrates the common logic of structure-driven thinking in number theory, pure mathematics, physics, and seismic imaging. The author also systematically proposes the Structure-Driven Inversion (SDI) paradigm, comprising Mathematical Structure-Driven Inversion (MSDI) and Physical Structure-Driven Inversion (PSDI). Structure-driven methodology does not negate objective-driven methodology but complements it: objective-driven is a general search strategy, whereas structure-driven is a custom direct strategy. The article concludes that the significance of this shift lies not only in efficiency gains but also in cognitive deepening---when facing complex problems, the first question is no longer ``how to solve the problem faster'' but ``what is the structure of the problem.''
\end{abstract}

\textbf{Keywords}: structure-driven; objective-driven; complex problems; Structure-Driven Inversion; mathematical structure; physical structure

\section{Starting from a Proof}

In September 2026, a formalized proof based on large-model automated reasoning was released on GitHub, proving a weakened form of the Liouville-variant Goldbach conjecture\cite{ref1}. The proof spans only two pages. Its core is not exhaustive enumeration or brute-force computation, but the identification of a sign structure of the Liouville function under additive composition, along which a concise proof path is directly constructed. The proof has been formally verified in Lean.

This event itself is a development at the intersection of number theory and artificial intelligence. Viewed against a broader background, however, it reveals a methodological trend emerging independently across multiple fields: \textbf{when problems become sufficiently complex, one rediscovers that identifying the intrinsic structure of a problem is more fundamental than blind search.}

The views expressed in this article stem from the author's forty years of learning, research, and application in inverse problems, spanning potential-field inversion, seismic imaging, and deep learning.

\section{The Two-Century Objective-Driven Tradition}

To understand this trend, one must first review the mainstream paradigm for solving inverse problems.

Since Laplace and Gauss, the dominant approach to solving inverse problems has been: \textbf{first construct an objective function, then solve by minimizing it.}

\begin{itemize}
    \item Laplace adopted the least-absolute-deviation criterion to fit astronomical observations\cite{ref2};
    \item Gauss proposed the least-squares method in 1809 to determine the orbit of Ceres\cite{ref3};
    \item Legendre independently proposed the least-squares method, published earlier than Gauss\cite{ref4};
    \item Tikhonov introduced regularization in the 1960s, stabilizing the solution of ill-posed problems\cite{ref5};
    \item The adjoint-state method extended it to large-scale geophysical inversion\cite{ref6,ref7};
    \item In recent years, SGD/Adam and other optimization algorithms reached new scales in deep learning\cite{ref8}.
\end{itemize}

Over more than two centuries, the form of the objective function and the minimization algorithms have changed, but the core logic of ``constructing an objective function and solving it'' has never changed.

This paradigm dominates inverse-problem solving because of its \textbf{universality}: as long as the forward operator is computable, whether differentiable or not, it can be incorporated into the framework of objective-function minimization.

Yet it also has three inherent limitations:

\begin{enumerate}
    \item \textbf{Computational efficiency bottleneck}---each iteration requires a complete forward computation, typically requiring tens to hundreds of iterations;
    \item \textbf{Ambiguous physical meaning}---the gradient is a mathematical construct; it reflects ``how the parameters should change,'' not the direct causal relationship between observed data and the physical process;
    \item \textbf{Multi-parameter mathematical coupling}---the cross-correlation operation mixes the contributions of different parameters, producing non-physical mathematical coupling that requires an expensive Hessian matrix for decoupling.
\end{enumerate}

These limitations are not accidental. Their root cause is that \textbf{the objective-driven paradigm treats all problems as the same ``differentiable function fitting'' problem, ignoring the intrinsic special structure of each problem.} When a translation-invariant operator is treated as a general matrix, and when a physical process with a causal chain is treated as a general mapping, efficiency, interpretability, and decoupling capability are all compromised.

This limitation has a deeper root. In the early twentieth century, Hadamard pointed out that the well-posedness of a mathematical problem depends on the existence, uniqueness, and continuous dependence of the solution on the data\cite{ref9}. Forward problems usually satisfy these three criteria; inverse problems usually do not. It was later recognized that inverse problems are not ``universally'' ill-posed---when a problem possesses special structure, the inversion can be well-posed\cite{ref10,ref11}. This is precisely the starting point of structure-driven thinking.

Notably, objective-driven thinking not only dominates inverse-problem solving in the natural sciences, but also pervades almost all domains of human practice, including economics, management, social governance, education, and healthcare: profit maximization, loss minimization, maximum probability, minimum risk, minimum error---all are different manifestations of objective-driven thinking. Their common feature is the definition of a quantifiable objective followed by some means of approaching the optimum. Yet \textbf{once the objective is set, the intrinsic structure of the system is often overlooked or even destroyed by objective-driven interventions.}

(Note: In the context of inverse-problem solving, objective-driven thinking typically manifests as optimization-driven thinking; the two are used interchangeably in this article. In discussing broader human practice, the term objective-driven is used.)

\section{Structure-Driven as a Complementary Paradigm}

In September 2026, the author published ``Structure-Driven Inversion: A New Paradigm for Solving Inverse Problems'' on arXiv\cite{ref12}, systematically proposing the SDI (Structure-Driven Inversion) paradigm.

The core idea of SDI is: \textbf{rather than approaching the solution through iterative search, identify and exploit the intrinsic structure of the problem to construct the solution method directly.}

It encompasses two types of structure:

\begin{itemize}
    \item \textbf{Mathematical structure}---the operator in the problem has a structure that can be unitarily diagonalized, such as translation invariance;
    \item \textbf{Physical structure}---the causal chain of the physical process in the problem has a structure that can be decomposed, such as wave propagation and the layered structure of deep networks.
\end{itemize}

And two corresponding driving inversion types:

\begin{itemize}
    \item \textbf{Mathematical Structure-Driven Inversion (MSDI)}: exploiting the unitary diagonalization property of the operator to construct a pseudo-inverse directly in the spectral domain\cite{ref13,ref14};
    \item \textbf{Physical Structure-Driven Inversion (PSDI)}: exploiting the causal chain of signal propagation to solve directly through virtual-source projection and deconvolution (inversion)\cite{ref15,ref16}.
\end{itemize}

Its core operation is \textbf{division/inversion/deconvolution}, rather than the objective-driven \textbf{cross-correlation/multiplication}. Division naturally decouples, avoiding the ``mathematical coupling'' introduced by multiplication in objective-driven methods, and eliminating the need for the Hessian matrix. SDI is the concrete practice of structure-driven methodology in inverse-problem solving.

The complete logical chain of structure-driven methodology can be summarized in four steps:

\begin{enumerate}
    \item \textbf{Deconstruct the problem}---decompose a complex problem into more basic components;
    \item \textbf{Identify the problem structure}---discover the intrinsic relationships among the decomposed components (orthogonality, causality, hierarchy, etc.);
    \item \textbf{Exploit the structure to construct a solution method}---construct a solution method directly based on the identified structure;
    \item \textbf{Solve the problem}---use the constructed solution method to solve the problem.
\end{enumerate}

Data processing based on harmonic analysis is a classic example of this logical chain: the Fourier transform performs the deconstruction, the independence of frequency components is the identified structure, the spectral filter is constructed, and one forward transform, one inverse transform, and spectral filtering complete the solution. More generally, wavelet decomposition, singular value decomposition, causal-chain decomposition, and hierarchical decomposition are all different ways of deconstructing a problem; they identify different types of structure, thereby constructing different solution methods. \textbf{The core of structure-driven thinking is to find the decomposition that matches the problem, and to identify and exploit the structure.}

It should be emphasized that structure-driven methodology is not a negation of objective-driven methodology. The two are complementary: \textbf{objective-driven is a ``general search strategy,'' while structure-driven is a ``custom direct strategy.''} When structure is missing, structure-driven degenerates to objective-driven; when structure is sufficient, structure-driven can provide concise and efficient solution methods.

\section{Cross-Domain Corroboration}

This ``structure-before-solution'' methodology is emerging independently across multiple fields.

\textbf{In number theory and artificial intelligence}, the aforementioned proof of the weakened Goldbach conjecture does not rely on exhaustive enumeration, but identifies the deep connection between the sign structure of the Liouville function and additive composition, and constructs a proof path along this structure\cite{ref1}. The proof was completed by a large-model automated reasoning system, demonstrating that structure identification is not only the intuition of human mathematicians, but can also be learned and executed by machines.

\textbf{In pure mathematics}, the proof of the three-dimensional Kakeya conjecture by Wang and Zahl is likewise a structure-driven example\cite{ref17}. Rather than directly addressing the original problem of ``the minimum area swept by a needle,'' they deconstruct it into the problem of ``the degree of overlap of thin tubes,'' identify structures such as multiscale nesting and convex-set density, and then use these structures to construct a multiscale inductive path, step by step pushing the dimension lower bound to 3. This work is methodologically identical to SDI: \textbf{deconstruct the problem, identify the structure, exploit the structure to construct a path.}

\textbf{In physics}, symmetry methods have always been a core tool. Noether's theorem reveals the profound connection between symmetry and conservation laws\cite{ref18}---physicists have always sought the symmetry structure of a problem and then used it to simplify the solution. From classical mechanics to quantum field theory, this logic runs throughout.

\textbf{In geophysical inversion}, structure-driven thinking runs through nearly four decades of the author's research, but potential-field inversion and seismic imaging rely on different types of structure.

\textbf{In potential-field inversion}, mathematical structure-driven thinking dominates. The potential-field forward operator is a translation-invariant convolutional operator that can be diagonalized by the Fourier transform. As early as the master's stage, the author discovered that the unitary structure of the Fourier transform matrix could be used directly to construct a pseudo-generalized inverse for convolutional linear inversion, published in 1993\cite{ref19}. Subsequently, the author exploited the same unitary diagonalization structure in gravity inversion, analytically extending 2D surface data to a 3D data volume, providing vertical wavenumber sampling for 3D gravity inversion, and constructing an efficient true 3D gravity inversion method\cite{ref14}.

\textbf{In seismic imaging}, physical structure-driven thinking dominates. Wave propagation is not a simple convolutional operator but a physical process with a causal chain: incident wave $\rightarrow$ virtual source $\rightarrow$ secondary wave. From constructing a local inverse operator directly using the locally space-invariant approximation in image-domain least-squares migration\cite{ref20}, to directly separating parameters using the physical law that reflection coefficients vary with angle in angle-domain relative impedance perturbation imaging\cite{ref21}, to waveform imaging based on linear forward representations of seismic data\cite{ref22,ref23} and isotropic elastic reverse time migration\cite{ref24}, physical structure-driven thinking runs throughout.

\textbf{In deep learning}, structure-driven thinking can likewise be extended to deep neural networks. PSDBP (Physical Structure-Driven Back-Propagation) analogizes a deep network to a physical system: inter-layer connections correspond to stiffness matrices, activation functions to constitutive relations, and forward propagation to an excitation-response process\cite{ref16}. Based on this analogy, PSDBP constructs three structural projections: inter-layer demand signals are transported through the transpose of the weights (geometric structure), the activation function is inverted analytically through $\sigma^{-1}$ to recover the target displacement (constitutive structure), and the stiffness matrix is updated through the input-residual correlation (least-squares projection). Its core operation is \textbf{function inversion}, rather than chain-rule differentiation. It should be noted that PSDBP does not aim to replace BP, but rather provides a physical-structure-based reinterpretation, showing that the structure-driven perspective and the chain-rule perspective can coexist and complement each other.

\textbf{In diffusion-model denoising}, structural priors combine with objective-driven thinking. In 2026, the author and his student, based on the observation that seismic signals have relatively simple structural characteristics, proposed changing the noise-fitting objective of the general diffusion model into a signal-fitting objective, and proposed the corresponding diffusion-model denoising method\cite{ref25}. In September of the same year, the author systematically distilled this practice, spanning nearly four decades, into the SDI paradigm\cite{ref12}. SDI originated during the author's master's studies, and the core has remained unchanged: \textbf{identify the structure, exploit the structure, and construct a solution method.}

The common logic across these fields is:

\begin{quote}
\textbf{Many problems already contain within themselves all the structural information needed to solve them. The real difficulty is not that ``computational power is insufficient,'' but that ``the structure of the problem has not been identified.'' Once the structure is correctly identified, the solution method will be concise and efficient.}
\end{quote}

\section{Methodological Significance}

From a methodological perspective, the emergence of structure-driven thinking reflects a cognitive shift:

\textbf{from ``how to search'' to ``how to identify.''}

Over the past half-century and more, the mainstream approach to inverse-problem solving has been: given a forward model and observed data, construct an objective function and approach the solution through iterative search. This is the objective-driven approach.

The structure-driven approach, by contrast, is: first identify the structure of the forward problem, then exploit the structure to construct the solution method for the inverse problem.

The significance of this shift lies not only in efficiency gains, but also in \textbf{cognitive deepening}: it prompts us to ask ``what is the intrinsic structure of the problem'' rather than ``how to search faster.''

Bayesian maximum a posteriori estimation is a typical representative of objective-driven thinking: it encodes structural information with a prior distribution $p(m)$, then combines it with the likelihood term into an objective function, solved by optimization. A Gaussian prior corresponds to Tikhonov regularization, a Laplace prior to L1 regularization, and a TV prior to total-variation regularization. However, \textbf{a prior distribution is only a ``soft'' statistical expression of structural information, not an explicit exploitation of the intrinsic structure of the operator or process.} The difference of structure-driven methodology is that it does not ``stuff'' structure into the objective function, but ``extracts'' structure to construct the solution method.

The claim of structure-driven methodology is: \textbf{before setting the objective, first ask what the intrinsic structure of the problem is; in pursuing the objective, exploit the structure rather than ignore it.} This principle applies equally to the natural sciences and the social sciences.

\section{Outlook}

Structure-driven methodology is not a closed theory but an open framework.

As new mathematical and physical structures are identified, more mathematics-driven and physics-driven methods can be developed. In number theory, pure mathematics, artificial intelligence, physics, and geophysical inversion, structure-driven thinking is also emerging in its own forms.

Perhaps we are witnessing the emergence of a more general methodology: \textbf{when facing complex problems, the first question is no longer ``how to compute faster,'' but ``what is the structure of the problem.''}

This shift may be more profound than any single algorithmic breakthrough.

The SDI paper is available at arXiv:2609.16536\cite{ref12}.

\end{document}